\documentclass[12pt,twoside]{amsart}
\usepackage{amsmath, amsthm, amssymb}

\date{}

\begin{document}

\centerline{}

\centerline {\Large{\bf Fuzzy Anti-norm and Fuzzy $\alpha$-anti-convergence}}

\centerline{}

%% My definition
\newcommand{\mvec}[1]{\mbox{\bfseries\itshape #1}}

\centerline{\bf {Bivas Dinda$^1$, T.K. Samanta$^2$ and Iqbal H. Jebril$^3$}}

\centerline{}

\centerline{$^1$Department of Mathematics,}
\centerline{ Mahishamuri Ramkrishna Vidyapith,}
\centerline{West Bengal, India. }
\centerline{E-mail: bvsdinda@gmail.com}
\centerline{}
\centerline{$^2$Department of Mathematics,}
\centerline{ Uluberia
College, West Bengal, India.}
\centerline{E-mail: mumpu$_{-}$tapas5@yahoo.co.in}
\centerline{}
\centerline{$^3$Department of Mathematics,}
\centerline{ King Faisal University, Saudi Arabia. }
\centerline{E-mail: iqbal501@hotmail.com}

\centerline{}

\newtheorem{Theorem}{\quad Theorem}[section]

\newtheorem{definition}[Theorem]{\quad Definition}

\newtheorem{theorem}[Theorem]{\quad Theorem}

\newtheorem{remark}[Theorem]{\quad Remark}

\newtheorem{corollary}[Theorem]{\quad Corollary}

\newtheorem{note}[Theorem]{\quad Note}

\newtheorem{lemma}[Theorem]{\quad Lemma}

\newtheorem{example}[Theorem]{\quad Example}

\centerline{{\bf Abstract}}
\textbf{\emph{In this paper the definition of fuzzy antinorm is modified. Some properties of finite dimensional fuzzy antinormed linear space are studied. Fuzzy $\alpha$-anti-convergence and fuzzy $\alpha$-anti-complete linear spaces are defined and some of their properties are studied. }}

%=============================
\section{Introduction}
%=============================
During the last few years there is a growing interest in the extension of fuzzy set theory which is a useful tool to describe the situation in which data are imprecise or vague or uncertain. Fuzzy set theory handle the situation by attributing a degree of membership to which a certain object belongs to a set. It has a wide range of application in the field of population dynamics \cite{Barros LC}, chaos control \cite{Fradkov}, computer programming \cite{Giles}, medicine \cite{Barro} etc. \\
The concept of fuzzy set theory was first introduce by Zadeh\cite{zadeh} in 1965
and thereafter, the concept of fuzzy set theory applied on different branches of pure and applied
mathematics in different ways by several authors. The concept of fuzzy norm was introduced by Katsaras
\cite{Katsaras} in 1984. In 1992, Felbin\cite{Felbin1} introduced
the idea of fuzzy norm on a linear space. Cheng-Moderson
\cite{Shih-chuan} introduced another idea of fuzzy norm on a linear
space whose associated metric is same as the associated metric of
Kramosil-Michalek \cite{Kramosil}. In 2003, Bag and Samanta
\cite{Bag} modified the definition of fuzzy norm of
Cheng-Moderson \cite{Shih-chuan} and established the concept of
continuity and boundednes of a linear operator with respect to
their fuzzy norm in \cite{Bag1}.\\
Later on  Jebril and Samanta \cite{Samanta1} introduced the concept of fuzzy anti-norm on a linear space depending on the idea of fuzzy anti norm was introduced by Bag and Samanta \cite{Bag2}. The motivation of introducing fuzzy anti-norm is to study fuzzy set theory with respect to the non-membership function. It is useful in the process of decision making.\\
In this paper we generalize the definition of fuzzy anti-norm on a linear space. Later on we prove Riesz lemma and some important properties of finite dimensional fuzzy anti-normed linear space.
Also, We define fuzzy $\alpha$-anti-convergence, fuzzy $\alpha$-anti-Cauchy sequence, fuzzy $\alpha$-anti-completeness and study the relation among them.

%==============================
\section{Preliminaries}
%==============================
This section contains some basic definition and preliminary results which will be needed in the sequel.
\begin{definition}
\cite{Schweizer,klement} A binary operation \, $\diamond \; : \; [\,0 \; ,
\; 1\,] \; \times \; [\,0 \; , \; 1\,] \;\, \rightarrow \;\,
[\,0 \; , \; 1\,]$ \, is continuous \, $t$-conorm if
\,$\diamond$\, satisfies the
following conditions \, $:$ \\
$(\,i\,)\;\;$ \hspace{0.1cm} $\diamond$ \, is commutative and
associative ,
\\ $(\,ii\,)\;$ \hspace{0.1cm} $\diamond$ \, is continuous , \\
$(\,iii\,)$ \hspace{0.1cm} $a \diamond 0 \,=\, a
\hspace{1.2cm}
\forall a  \in [\,0 , 1\,]$ , \\
$(\,iv\,)$ \hspace{0.1cm} $a \diamond b \, \leq \, c
\diamond d$ \, whenever \, $a \leq c$  ,  $b \leq d$
 and  $a,\,b,\,c,\,d\in [\,0
, 1\,].$
\end{definition}
A few examples of continuous t-conorm are $\,a\,\diamond\,b\,=\,a+b-ab,\;\,a\,\diamond\,b\,=\,\max\{a,b\},\;\,a\,\diamond\,b\,=\,\min\{a+b,1\}$.

\begin{definition}\label{p1}
\cite{Bag2} Let $X$ be a linear space over $F$ (field of real/complex numbers). Let $N^*$ be a fuzzy subset of $X\times \mathbb{R}$ such that for all $x,\,y\in X$ and $c\;\in\;F$\\
$(N^*1)$ $\;\forall t\,\in \,\mathbb{R}$ with $t\,\leq\,0\;,\;N^*\,(\,x\,,\,t\,)\;=\;1\;$,\\
$(N^*2)$ $\;\forall t\,\in \,\mathbb{R}$ with $t\,>\,0\;,\;N^*\,(\,x\,,\,t\,)\;=\;0\;$ if and only if $x\;=\;\theta$,\\
$(N^*3)$ $\;\forall t\,\in \,\mathbb{R}$ with $t\,>\,0,\,N^*\,(\,cx\,,\,t\,)=\,N^*\,(\,x\,,\,\frac{t}{\mid c \mid}\,)$ if $c\neq\;0\,$,\\
$(N^*4)$ $\;\forall s,t\,\in \,\mathbb{R}$ with $\;N^*(x\,+\,y\,,\,s\,+\,t)\;\leq\; \max \;\{\;N^*\,(\,x\,,\,s\,)\;,\;N^*\,(\,y\,,\,t\,)\;\}$,\\
$(N^*5)$ $N^*(\,x\,,\,\cdot\,)$ is a non-increasing of $t\in \mathbb{R}$ and $\mathop{\lim }\limits_{t\to\infty}N^*(x,t)=0.$
\\Then $\,N^*$ is called a B-S-fuzzy antinorm on $X$.
\end{definition}

We assume that\\
$(N^*6)$ For all $t\in\mathbb{R}$ with $t\,>\,0\;,\;N^*\,(\,x\,,\,t\,)\;<\;1\;$ implies $x\;=\;\theta$.

\begin{definition}
\cite{Samanta1} Let, $(\,U\,,\,N^*\,)$ be a B-S-fuzzy antinormed linear space. A sequence $\{x_n\}_{n\in\mathbb{N}}$ in $U$ is said to \textbf{converge} to $x\,\in\,U$ if given $t\,>\,0\;,\;r\,\in\,(\,0\,,\,1\,)$ there exist an integer $n_0\,\in\,\mathbb{N}$ such that
\[N^*\,(\;x_n\,-\,x\,,\,t\,)\;<\;r\;\;\;\forall\;n\,\geq\,n_0.\]
\end{definition}

\begin{definition}
\cite{Samanta1} Let, $(\,U\,,\,N^*\,)$ be a B-S-fuzzy antinormed linear space. A sequence $\{x_n\}_{n\in\mathbb{N}}$ in $U$ is said to be \textbf{Cauchy sequence} if for given $t\,>\,0\;,\;r\,\in\,(\,0\,,\,1\,)$ there exists an integer $n_0\,\in\,\mathbb{N}$ such that
\[N^*\,(\;x_{n+p}\,-\,x_n\,,\,t\,)\;<\;r\;\;\;\forall\;n\,\geq\,n_0\;\,,\,p=1,2,3,\cdots.\]
\end{definition}

\begin{definition}
\cite{Samanta1} A subset $A$ of a B-S-fuzzy antinormed linear space $(\,U\,,\,N^*\,)$ is said to be \textbf{bounded} if and only if there exist $t\,>\,0\;,\;r\,\in\,(\,0\,,\,1\,)$ such that \[N^*\,(\,x\,,\,t\,)\;<\;r\;\;\;\forall\;x\,\in\,A.\]
\end{definition}

\begin{definition}
\cite{Samanta1} A subset $A$ of a B-S-fuzzy antinormed linear space $(\,U\,,\,N^*\,)$ is said to be \textbf{compact} if any sequence $\{x_n\}_{n\in\mathbb{N}}$ in $A$ has a subsequence converging to an element of $A$.
\end{definition}

\begin{definition}
\cite{Samanta1} Let $(\,U\,,\,N^*\,)$ be a B-S-fuzzy antinormed linear space. A subset $B$ of $U$ is said to be \textbf{closed} if any sequence $\{x_n\}_{n\in\mathbb{N}}$ in $B$ converges to $x\,\in\,B$, that is, \[\mathop {\lim }\limits_{n\; \to \;\infty }\,N^*\,(\;x_n\,-\,x\,,\,t\,)\;=\;0,\;\;\forall\,t\,>\,0\;\Rightarrow\;\,x\,\in\,B.\]
\end{definition}

%==========================================
\section{Fuzzy Anti-normed Linear Space}
%==========================================
In this section the definition of B-S-fuzzy antinorm is modified and after modification it will be termed as fuzzy antinorm with respect to a t-conorm $\diamond$. Thereafter some important results will be deduced.
\begin{definition}\label{d1}
Let $V$ be linear space over the field $F\,(=\,\mathbb{R}\;or\;\mathbb{C}\,)$. A fuzzy subset $\;\nu\;$ of $V\,\times\,\mathbb{R}$ is called a \textbf{fuzzy antinorm} on $V$ with respect to a t-conorm $\diamond$ if and only if for all $x\,,\,y\;\in\;V$ \\
$(i)$ $\;\;\forall t\,\in \,\mathbb{R}$ with $t\,\leq\,0\;,\;\nu(\,x\,,\,t\,)\;=\;1\;$;\\
$(ii)$ $\;\,\forall t\,\in \,\mathbb{R}$ with $t\,>\,0\;,\;\nu(\,x\,,\,t\,)\;=\;0\;$ if and only if $x\;=\;\theta$;\\
$(iii)$ $\,\forall t\,\in \,\mathbb{R}$ with $t\,>\,0\;,\;\nu\,(\,cx\,,\,t\,)=\,\nu\,(\,x\,,\,\frac{t}{\mid c \mid}\,)\;$ if $c\neq\;0\,,c\;\in F$;\\
$(iv)$ $\;\forall s,\,t\;\in\;\mathbb{R}$ with $\;\nu\,(\,x\,+\,y\,,\,s\,+\,t\,)\;\leq\;\nu\,(\,x\,,\,s\,)\;\diamond\;\nu\,(\,y\,,\,t\,)\;$;\\
$(v)$  $\mathop{\lim }\limits_{t\, \to\,\;\infty }\,\nu\,(\,x\,,\,t\,)=0.$
\end{definition}

The definition \ref{d1} is more general than the definition \ref{p1}; since, in $(N^*4)$ instead of maximum function we have used more generalized function, conorm function and in the condition $(N^*5)$ it is used that $N^*(x,\cdot)$ is non-increasing function of $t(\in \mathbb{R}),$ which is redundant and later on it will be deduced.

\begin{remark}\label{r1}
Let  $\nu$  be a fuzzy anti-norm on $V$ with respect to a t-conorm $\diamond$ then $\nu(x,t)$ is non-increasing with respect to $t$ for each $x\,\in\,V$.
\end{remark}
{\bf Proof.}$\;\;$
Let $t\,<\,s$. Then $k\,=\,s-t\,>\,0$ we have
\[\nu(x,t)\,=\,\nu(x,t)\,\diamond\,0\,=\,\nu(x,t)\,\diamond\,\nu(0,k)\,\geq\,\nu(x,s).\]
Hence the proof.

\begin{definition}
If $\;A^*\;=\;\{\,((\,x\,,\,t\,)\;,\;\nu\,(\,x\,,\,t\,))\;:\;(\,x\,,\,t\,)\;\in\;V\,\times\,\mathbb{R}\}$ is a fuzzy antinorm on a linear space $V$ with respect to a t-conorm $\diamond$ over a field $F$, then $(\,V\,,\,A^*\,)$ is called a fuzzy antinormed linear space with respect to the t-conorm $\diamond$ over the field $F$.
\end{definition}

We further assume that for any fuzzy anti-normed linear space $\,(\,V\,,\,A^*\,)$ with respect to a t-conorm $\diamond$,\\
{\bf (vi)}$\;\;\,\nu\,(\,x\,,\,t\,)\,<\,1\;,\;\forall\,t>0\;\Rightarrow\;x\,=\,\theta.\\ $
{\bf (vii)}$\;\,\nu(\,x\,,\,\cdot\,)\,$ is a continuous function of $\mathbb{R}$ and strictly decreasing on the subset $\,\{\,t\,:\;0<\nu(x,t)<1\}\;$ of $\,\mathbb{R}.\\$
{\bf (viii)} $a\diamond a=a,\;\;\forall a\in[0,1].$

\begin{example}
Let $(V,\|\cdot\|)$ be a normed linear space and consider $a\diamond b=a+b-ab$. Define $\nu:V\times \mathbb{R}\rightarrow [0,1]\,$ by
\[\nu\,(\,x\,,\,t\,)\;=\begin{cases}\;0,\;\;\;\,if\;t\,>\,\|\,x\,\|\\
1,\;\;\;\,if\,\;t\,\leq\,\|\,x\,\|.\;\hspace{-1.8 cm}
\end{cases}\]
Then $\nu$ is a fuzzy antinorm on $V$ with respect to the t-conorm $\diamond$ and $(\;V\;,\;\nu\;)$ is a fuzzy anti-normed linear space with respect to the t-conorm $\diamond$.
\end{example}
{\bf Solution.} $(i)\;\;\forall x\in V$ and $\forall t\in \mathbb{R},\,t\leq 0$ we have $\nu(x,t)=1.$\\
$(ii)\;\;\forall t\in \mathbb{R},\,t>0\,$ we have $\nu(\theta,t)=0.$ Again
\[\nu(x,t)=0,\,\forall t>0\,\Leftrightarrow\,\|x\|<t,\,\forall t(>0)\in \mathbb{R}\,\Leftrightarrow\,\|x\|=0\,\Leftrightarrow\,x=\theta.\]
\[(iii)\;\;\nu(cx,t)=0\,\Leftrightarrow\,t>\|cx\|\,\Leftrightarrow\,t>|c|\|x\|\,\Leftrightarrow\,\frac{t}{|c|}>\|x\|\, \Leftrightarrow\,\nu(x,\frac{t}{|c|})=0.\]
\[\nu(cx,t)=1\,\Leftrightarrow\,t\leq\|cx\|\,\Leftrightarrow\,t\leq|c|\|x\|\,\Leftrightarrow\,\frac{t}{|c|}\leq\|x\|\, \Leftrightarrow\,\nu(x,\frac{t}{|c|})=1.\]
$(iv) \;\;\nu(x,s)\diamond\nu(y,t)=\nu(x,s)+\nu(y,t)-\nu(x,s)\nu(y,t).\\$
If $s>\|x\|$ and $t>\|y\|$ then
$\nu(x+y,s+t)=0\,,\;\;\;$ since $s+t>\|x\|+\|y\|\\$
and $\nu(x,s)\diamond \nu(y,t)=0+0-0=0.\\$
So, $\nu(x+y,s+t)=\nu(x,s)\diamond \nu(y,t).\\$
If $s>\|x\|$ and $t\leq\|y\|$ then $\nu(x,s)\diamond \nu(y,t)=0+1-0=1.\\$
If $s\leq\|x\|$ and $t>\|y\|$ then $\nu(x,s)\diamond \nu(y,t)=1+0-0=1.\\$
If $s\leq\|x\|$ and $t\leq\|y\|$ then $\nu(x,s)\diamond \nu(y,t)=1+1-1=1.\\$
Therefore in any of the above three cases $\nu(x,s)\diamond \nu(y,t)=1\geq \nu(x+y,s+t).\\$
Thus \[\nu(x+y,s+t)\leq \nu(x,s)\diamond \nu(y,t).\]
$(v)$  From the definition it is clear that $\mathop{\lim }\limits_{t\, \to\,\;\infty }\,\nu\,(\,x\,,\,t\,)=0.\\$
Thus $\nu$ is a fuzzy antinorm on $V$ with respect to the t-conorm $\diamond$ and $(\;V\;,\;\nu\;)$ is a fuzzy anti-normed linear space with respect to the t-conorm $\diamond$.

\begin{note}
The above example satisfes the condition (vi) but does not satisfy the condition (vii).
\end{note}

\begin{example}
Let $(V,\|\cdot\|)$ be a normed linear space and consider $a\diamond b=\min\{a+b,1\}$. Define $\nu:V\times \mathbb{R}\rightarrow [0,1]\,$ by
\[\nu\,(\,x\,,\,t\,)\;=\begin{cases}0,\;\;\;\,\hspace{2 cm}if\;t\,>\,\|\,x\,\|,\,t>0.\hspace{0 cm}
\\ \frac{\|x\|}{t+\|x\|},\;\;\;\hspace{1.3 cm}if\,\;t\,\leq\,\|\,x\|\,,\,t>0.\;\hspace{-1.8 cm}
\\1,\;\;\;\hspace{2 cm}if\,\;t\,\leq\,0.\;\hspace{0 cm}\end{cases}\]
Then $\nu$ is a fuzzy antinorm on $V$ with respect to the t-conorm $\diamond$ and $(\;V\;,\;\nu\;)$ is a fuzzy anti-normed linear space with respect to the t-conorm $\diamond$.
\end{example}
{\bf Solution.} $(i)$  From the definition we have $\nu(x,t)=1$ if $t\leq 0,\,\forall t\in \mathbb{R}.\\$
$(ii)$  If $t>0$ and $t>\|x\|$ then
\[\nu(x,t)=0\,\Leftrightarrow\,\|x\|<t,\,\forall t(>0)\in \mathbb{R}\,\Leftrightarrow\,\|x\|=0\,\Leftrightarrow\,x=\theta.\]
If $t\leq 0$ and $t\leq\|x\|$ then
\[\nu(x,t)=0\,\Leftrightarrow\,\frac{\|x\|}{t+\|x\|}=0\,\Leftrightarrow\,\|x\|=0\,\Leftrightarrow\,x=\theta.\]
\[(iii)\;\;\nu(cx,t)=0\,\Leftrightarrow\,t>\|cx\|\,\Leftrightarrow\,t>|c|\|x\|\,\Leftrightarrow\,\frac{t}{|c|}>\|x\|\, \Leftrightarrow\,\nu(x,\frac{t}{|c|})=0.\]
\[\nu(cx,t)=\frac{\|cx\|}{t+\|cx\|}\,\Leftrightarrow\,t\leq \|cx\|\,\Leftrightarrow\,\frac{t}{|c|}\leq \|x\|\,\Leftrightarrow\,\nu(x,\frac{t}{|c|}) =\frac{\|x\|}{\frac{t}{|c|}+\|c\|}=\frac{\|cx\|}{t+\|cx\|}.\]
\[(iv)\;\;\;\nu(x,s)\diamond\nu(y,t)=\min\{\nu(x,s)+\nu(y,t),1\}.\hspace{5 cm}\]
If $\|x\|\geq s$ and $\|y\|\geq t$ then
\[\nu(x,s)+\nu(y,t)=\frac{\|x\|}{s+\|x\|}\,+\,\frac{\|y\|}{t+\|y\|}\hspace{9 cm}\]
\[=\frac{(t\|x\|+\|x\|\|y\|+s\|y\|)+\|x\|\|y\|}{(t\|x\|+\|x\|\|y\|+s\|y\|)+st}\geq \,1\;\;since, \;\|x\|\|y\|\geq st.\]
In this case  $\nu(x,s)\diamond\nu(y,t)=1\,\geq\, \nu(x+y,s+t).\\$
If $\|x\|\geq s$ and $\|y\|<t$ then either $\|x+y\|\geq s+t$  or  $\|x+y\|<s+t.\\$
Now, $\nu(x,s)+\nu(y,t)=\frac{\|x\|}{s+\|x\|}\,+0\,<\,1.\;$ Hence $\nu(x,s)\diamond\nu(y,t)=\frac{\|x\|}{s+\|x\|}.\\$
If $\|x+y\|\geq s+t$  then
\[\nu(x+y,s+t)-\nu(x,s)\,\diamond\,\nu(y,t)=\frac{\|x+y\|}{s+t+\|x+y\|}-\frac{\|x\|}{s+\|x\|}\hspace{4 cm}\]
\[\leq \frac{\|x\|+\|y\|}{s+t+\|x\|+\|y\|}-\frac{\|x\|}{s+\|x\|}=\frac{s\|y\|-t\|x\|}{(s+t+\|x\|+\|y\|)(s+\|x\|)}\hspace{2 cm}\]
\[<\frac{st-t\|x\|}{(s+t+\|x\|+\|y\|)(s+\|x\|)},\;\;since\;\|y\|<t.\hspace{5 cm}\]
\[\leq\,0,\hspace{5 cm}since\;s\leq\, \|x\|\Rightarrow\,st<t\|x\|.\hspace{6 cm}\]
Therefore, $\nu(x+y,s+t)<\nu(x,s)\diamond\nu(y,t).\\$
If $\|x+y\|<s+t$ then
\[\nu(x+y,s+t)=0\leq\frac{\|x\|}{s+\|x\|}=\nu(x,s)\diamond\nu(y,t).\]
If $\|x\|< s$ and $\|y\|\geq t$ then in the similar manner (as in the case when $\|x\|\geq s$ and $\|y\|<t$) we can show that $\nu(x+y,s+t)\leq\nu(x,s)\diamond\nu(y,t).$
If $\|x\|< s$ and $\|y\|<t$ then $\nu(x,s)+\nu(y,t)=0+0<1.\;$ Therefore, $\nu(x,s)\diamond\nu(y,t)=0.\\$
Also $\|x+y\|\leq\|x\|+\|y\|<s+t\;$ and hence $\nu(x+y,s+t)=0.\\$
Hence $\;\nu(x+y,s+t)=\nu(x,s)\diamond\nu(y,t).\\$
Thus  in any case \[\nu(x+y,s+t)\leq \nu(x,s)\diamond\nu(y,t)\]
(v)  If $t>\|x\|$ then from the definition it is clear that $\mathop{\lim }\limits_{t\, \to\,\;\infty }\,\nu\,(\,x\,,\,t\,)=0.\\$
If $x\neq\theta$ and $t\leq\|x\|$ then \[\mathop{\lim }\limits_{t\, \to\,\;\infty }\,\nu\,(\,x\,,\,t\,)=\mathop{\lim }\limits_{t\, \to\,\;\infty }\frac{\|x\|}{t+\|x\|}=0.\]
If $x=\theta$ and $t\leq\|x\|$ then \[\mathop{\lim }\limits_{t\, \to\,\;\infty }\,\nu\,(\,x\,,\,t\,)=\mathop{\lim }\limits_{t\, \to\,\;\infty }\nu(\theta,t)=\mathop{\lim }\limits_{t\, \to\,\;\infty }\frac{1}{t}=0.\]
Hence $\mathop{\lim }\limits_{t\, \to\,\;\infty }\,\nu\,(\,x\,,\,t\,)=0\;\;\;\forall x\in V.\\$
Thus $\nu$ is a fuzzy antinorm on $V$ with respect to the t-conorm $\diamond$ and $(\;V\;,\;\nu\;)$ is a fuzzy anti-normed linear space with respect to the t-conorm $\diamond$.

\begin{note}
The above example does not satisfy the conditions (vi) and (vii).
\end{note}

\begin{example}
Let $(V,\|\cdot\|)$ be a normed linear space and consider $a\diamond b=\max\{a,b\}$. Define $\nu:V\times \mathbb{R}\rightarrow [0,1]\,$ by
\[\nu\,(\,x\,,\,t\,)\;=\begin{cases}\frac{\|x\|}{t+\|x\|},\;\;\;\hspace{1.3 cm}if\,\;\,t>0.\;\hspace{-1.8 cm}
\\1,\;\;\;\hspace{2 cm}if\,\;t\,\leq\,0.\;\hspace{0 cm}\end{cases}\]
Then by Example 3.2 of \cite{Samanta1} it follows that $\nu$ is a fuzzy antinorm on $V$ with respect to the t-conorm $\diamond$ and $(\;V\;,\;\nu\;)$ is a fuzzy anti-normed linear space with respect to the t-conorm $\diamond$.
\end{example}

\begin{note}
The above example does not satisfy the condition (vi) and satisfies the condition (vii).
\end{note}

\begin{example}
Let $(V,\|\cdot\|)$ be a normed linear space and consider $a\diamond b=\min\{\,a+b\,,\,1\,\}$. Define $\nu:V\times \mathbb{R}\rightarrow [0,1]\,$ by
\[\nu\,(\,x\,,\,t\,)\;=\begin{cases}\frac{\|x\|}{\;\,2t-\|x\|\,},\;\;\;\hspace{1 cm}if\,\;\,t>\|x\|.\;\hspace{-1.8 cm}
\\1,\;\;\;\hspace{2 cm}if\,\;t\,\leq\,\|x\|.\;\hspace{0 cm}\end{cases}\]
Then $\nu$ satisfies all the condition of Definition \ref{d1}. Therefore $\nu$ is a fuzzy antinorm on $V$ with respect to the t-conorm $\diamond$ and $(\;V\;,\;\nu\;)$ is a fuzzy anti-normed linear space with respect to the t-conorm $\diamond$.
\end{example}

\begin{note}
The above example satisfies both the conditions (vi) and (vii).
\end{note}

\begin{theorem}
Let $(\,V\,,\,A^*\,)$ be a fuzzy antinormed linear space with respect to a t-conorm $\diamond$ satisfying $(vi)$ and $(viii)$. Then for any $\alpha\;\in\;(\,0\;,\;1\,)$ the function $\left\| {\;x\;} \right\|_{\,\alpha }^{\,\ast}:X\rightarrow [0,\infty)$ defined as \[{\bf (ix)}\hspace{2.5 cm}\left\| {\;x\;} \right\|_{\,\alpha }^{\,\ast}\;=\;\bigwedge\,\{\,t>0\,:\,\nu\,(\,x,t\,)\leq 1\,-\,\alpha\},\;\alpha\in(\,0,1\,)\] is a norm on $V$.
\end{theorem}
{\bf Proof.}(i) For $x\in V,\\$
$\nu(x,t)=1\;$ for $t\leq 0\;\Rightarrow\,\bigwedge\,\{\,t>0\,:\,\nu\,(\,x,t\,)\leq 1\,-\,\alpha\}\,\geq\,0,\;\alpha\in(\,0,1\,)\;\\\Rightarrow\,\left\| {\;x\;} \right\|_{\,\alpha }^{\,\ast}\geq 0,\;\alpha\in(\,0,1\,).\\\\$
(ii) $\left\| {\;x\;} \right\|_{\,\alpha }^{\,\ast}=\, 0\;\Rightarrow\,\nu\,(\,x,t\,)\leq 1\,-\,\alpha\,<\,1,\;\forall\,t\in\mathbb{R},\,t>0\;\\\Rightarrow\,x=\theta,\;$ [by (vi)].\\
Conversely, $x=\theta\;\Rightarrow\,\nu(x,t)=0,\;\forall\,t>0\;\Rightarrow\,\bigwedge\, \{\,t>0\,:\,\nu\,(\,x,t\,)\leq 1\,-\,\alpha\}\,=\,0,\;\forall\,\alpha\in(\,0,1\,)\;\Rightarrow\,\left\| {\;x\;} \right\|_{\,\alpha }^{\,\ast}=0.\\\\$
(iii) If $c\neq 0$ then
\[\left\| {\;cx\;} \right\|_{\,\alpha }^{\,\ast}\;=\;\bigwedge\,\{\,s>0\,:\,\nu\,(\,cx,s\,)\leq 1\,-\,\alpha\}\,\hspace{1.5 cm}\]
\[=\;\bigwedge\,\{\,s>0\,:\,\nu\,(\,x,\frac{s}{|c|}\,)\leq 1\,-\,\alpha\}\,\hspace{0 cm}\]
\[= \;\bigwedge\,\{\,|c|t>0\,:\,\nu\,(\,x,t\,)\leq 1\,-\,\alpha\}\,\hspace{0 cm}\]
\[=\;\bigwedge\,|c|\,\{\,t>0\,:\,\nu\,(\,x,t\,)\leq 1\,-\,\alpha\}\,\hspace{0 cm}\]
\[=\;|\,c\,|\,\left\| {\;x\;} \right\|_{\,\alpha }^{\,\ast}.\hspace{4 cm}\]
If $c=0$ then $\left\| {\;cx\;} \right\|_{\,\alpha }^{\,\ast}\,=\,\left\| {\;\theta\;} \right\|_{\,\alpha }^{\,\ast}\,=\,0\,=\,0.\left\| {\;x\;} \right\|_{\,\alpha }^{\,\ast}\,=\;|\,c\,|\,\left\| {\;x\;} \right\|_{\,\alpha }^{\,\ast}.\\\\$
(iv) $\left\| {\;x\;} \right\|_{\,\alpha }^{\,\ast}\,+\,\left\| {\;y\;} \right\|_{\,\alpha }^{\,\ast}\,$\[=\;\bigwedge\,\{\,t>0\,:\,\nu\,(\,x,t\,)\leq 1\,-\,\alpha\}\,+\,\;\bigwedge\,\{\,s>0\,:\,\nu\,(\,y,s\,)\leq 1\,-\,\alpha\}\;\forall\,\alpha\in(0,1)\]
\[\geq\;\bigwedge\,\{\,t+s>0\,:\,\nu\,(\,x,t\,)\leq 1\,-\,\alpha,\;\nu\,(\,y,s\,)\leq 1\,-\,\alpha\}\hspace{4 cm}\]
\[\geq\;\bigwedge\,\{\,t+s>0\,:\,\nu\,(\,x+y,t+s\,)\leq 1\,-\,\alpha\;\;[by\;(viii)]\hspace{4 cm}\]
\[=\left\| {\;x+y\;} \right\|_{\,\alpha }^{\,\ast}.\hspace{12 cm}\]
Hence $\{\left\| {\;\cdot\;} \right\|_{\,\alpha }^{\,\ast}\}$ is a norm on $V.$

\begin{remark}
The norm defined above is more general than the norm defined in Theorem 3.2 0f \cite{Bag2}; since instead of $\nu(x,t)<\alpha$ we write $\nu(x,t)\leq 1-\alpha$.
\end{remark}

\begin{theorem}
Let $(\,V\,,\,A^*\,)$ be a fuzzy antinormed linear space with respect to a t-conorm $\diamond$. If $\alpha_1\leq \alpha_2$, then $\left\| {\;x\;} \right\|_{\,\alpha_1 }^{\,\ast}\leq \left\| {\;x\;} \right\|_{\,\alpha_2 }^{\,\ast}\; i.e.,\;\{\left\| {\;\cdot\;} \right\|_{\,\alpha }^{\,\ast}:\alpha\in(0,1)\}$ is a increasing family of norm on $V$.
\end{theorem}
{\bf Proof.} $\alpha_1\leq \alpha_2$ we have\\
\[\{\,t>0\,:\,\nu\,(\,x,t\,)\leq 1\,-\,\alpha_2\}\,\subset\,\{\,t>0\,:\,\nu\,(\,x,t\,)\leq 1\,-\,\alpha_1\}\hspace{1 cm}\]
\[\Rightarrow\;\bigwedge\{\,t>0\,:\,\nu\,(\,x,t\,)\leq 1\,-\,\alpha_2\}\,\geq\,\bigwedge\{\,t>0\,:\,\nu\,(\,x,t\,)\,\leq\, 1\,-\,\alpha_1\}\hspace{5 cm}\]
\[\Rightarrow\;\left\| {\;x\;} \right\|_{\,\alpha_2}^{\,\ast}\,\geq\,\left\| {\;x\;} \right\|_{\,\alpha_1 }^{\,\ast}\hspace{10 cm}\]

In the following theorem we describe another one equivalent expression for $\nu$, which will be useful to describe Riesz theorem in fuzzy environment.

\begin{theorem}\label{th1}
Let, $(\,V\,,\,A^*\,)$ be a fuzzy antinormed linear space with respect to a t-conorm $\diamond$ satisfying $(vi),\,(vii),\,(viii)$ and let $\;\nu\,^\prime\,:\,V\,\times\,R\;\rightarrow\;[\,0\,,\,1\,]$ be defined by \[{\bf (x)}\hspace{0.25 cm}\nu\;^\prime\,(\,x\,,\,t\,)\;=\begin{cases}\;\bigwedge\,\{\,1\,-\,\alpha\;:\;\left\| {\;x\;} \right\|_{\,\alpha }^{\,\ast}\;\leq\;t\},\hspace{0.25 cm}if\;(\,x\,,\,t\,)\neq(\,\theta\,,\,0\,)\\
\;1,\hspace{4.5 cm}\;if\;(\,x\,,\,t\,)=(\,\theta\,,\,0\,)\hspace{-2 cm}\end{cases}\]
Then $\nu\,^\prime\;=\;\nu,\;$ where $\left\| {\;\cdot\;} \right\|_{\,\alpha }^{\,\ast}$ is a increasing family of norms given by $(ix)$.
\end{theorem}
To prove this theorem we use the following lemma.
\begin{lemma}\label{l1}
Let  $(\,V\,,\,A^*\,)$ be a fuzzy antinormed linear space with respect to a t-conorm $\diamond$ satisfying $(vi),\,(vii),\,(viii)$ and $\{\;\left\| {\;\cdot\;} \right\|_{\,\alpha }^{\,\ast}\;:\;\alpha\;\in\;(\,0\;,\;1\,)\;\}$ be increasing family of norms of V, defined by $(ix)$. Then for $x_0(\;\neq\;\theta)\;\in\;V\;,\;\alpha\;\in\;(\,0\;,\;1\,)$ and $s(>0)\in\mathbb{R},$
\[\left\| {\;x_0\;} \right\|_{\,\alpha }^{\,\ast}\;=\;s\;\Leftrightarrow\;\nu\,(\,x_0\,,\,s\,)=1\,-\,\alpha.\;\]
\end{lemma}
{\bf Proof.}$\;\;$
Let $\left\| {\;x_0\;} \right\|_{\,\alpha }^{\,\ast}\;=\;s$, then $s\,>\,0$. Then there exist a sequence $\{s_n\}_{n\in\mathbb{N}}$, $s_n\,>\,0$ such that $\nu\,(\;x_0\,,\,s_n\;)\;\leq\;1\,-\,\alpha\;$, for all $n\,\in\,\mathbb{N}$ and $s_n\;\rightarrow\;s$ as $n\;\rightarrow\;\infty$. Therefore
\[\mathop {\lim }\limits_{n\; \to \;\infty }\,\nu\,(\;x_0\,,\,s_n\;)\;\leq\;1\,-\,\alpha\;\Rightarrow\;\nu\,(\;x_0\,,\,\mathop {\lim }\limits_{n\; \to \;\infty }\,s_n\;)\;\leq\;1\,-\,\alpha\;\;\;\;by\;(vii)\]
\[\Rightarrow\;\nu\,(\,x_0\;,\;\left\| {\;x_0\;} \right\|_{\,\alpha }^{\,\ast}\,)\;\leq\;1\,-\,\alpha\;,\;\forall\;\alpha\in(\,0,1\,).\hspace{-5.6 cm}\]
Let $\alpha\;\in\;(\,0\;,\;1\,)\;,\;x_0(\neq\,\theta\;)\in\;V$ and $s\;=\;\left\| {\;x_0\;} \right\|_{\,\alpha }^{\,\ast}\;=\;\bigwedge\,\{\,t\;:\;\nu\,(\,x_0\;,\;t\,)\;\leq\;1\,-\,\alpha\}\,.\;$ Since $\nu(x,\cdot)$ is continuous (by (vii)) we have
\begin{equation}\label{eq1}
\nu\,(\,x_0\,,\,s\,)\;\leq\;1\,-\,\alpha.
\end{equation}
If possible let $\nu\,(\,x_0\,,\,s\,)\;<\;1\,-\,\alpha$, then by (vii) there exist $s\,^\prime>s$ such that $\nu\,(\,x_0\,,\,s\,^\prime\,)<\nu\,(\,x_0\,,\,s\,)<\;1\,-\,\alpha$, which is impossible since $s\;=\;\bigwedge\,\{\,t\;:\;\nu\,(\,x_0\;,\;t\,)\;\leq\;1\,-\,\alpha\}\;$. Thus
\begin{equation}\label{eq2}
\nu\,(\,x_0\,,\,s\,)\;\geq\;1\,-\,\alpha.
\end{equation}
From (\ref{eq1}) and (\ref{eq2}) it follows that $\nu\,(\,x_0\,,\,s\,)\;=\;1\,-\,\alpha.\;\,$Thus
\begin{equation}\label{eq3}
\left\| {\;x_0\;} \right\|_{\,\alpha }^{\,\ast}\;=\;s\;\;\Rightarrow\;\nu\,(\,x_0\,,\,s\,)\;=\;1\,-\,\alpha.
\end{equation}
Next if $\nu\,(\,x_0\,,\,s\,)\;=\;1\,-\,\alpha\;,\;\alpha\;\in\;(\,0\;,\;1\,)$ then by $(vii)$
\begin{equation}\label{eq4}
\left\| {\;x_0\;} \right\|_{\,\alpha }^{\,\ast}\;=\;\bigwedge\,\{\,t\;:\;\nu\,(\,x_0\;,\;t\,)\;\leq\;1\,-\,\alpha\}\;=\;s.
\end{equation}
Hence from (\ref{eq3}) and (\ref{eq4}) we have for $\alpha\;\in\;(\,0\;,\;1\,)\;,\;x\,(\neq\,\theta\;)\;\in\;V$ and for $s\,>\,0\;,\;\;\left\| {\;x_0\;} \right\|_{\,\alpha }^{\,\ast}\;=\;s\;\Leftrightarrow\;\nu\,(\,x_0\,,\,s\,)\;=\;1\,-\,\alpha.\\\\$
\textit{Proof of the main theorem:\\}
Let $(x_0,t_0)\in V\times\mathbb{R}$. To prove this theorem we consider the following cases:\\
\textbf{case-1:$\;$}For any $\;x_0\,\in V\;$ and $\;t\leq 0,\,\;\nu\,(\,x_0\,,\,t_0\,)\;=\;\nu\,^\prime(\,x_0\,,\,t_0\,)\;=\;1.$\\
\textbf{case-2:$\;$}If $\;x_0\;=\;\theta\;,t_0\;>\;0\;.\;$ Then $\;\nu\,(\,x_0\,,\,t_0\,)\;=\;\nu\,^\prime(\,x_0\,,\,t_0\,)\;=\;0.$\\
\textbf{case-3:$\;$}$\;x_0\;\neq\;\theta\;,t_0\;(\,>\;0\,)\;\in\;\mathbb{R}$  such that $\nu\,(\,x_0\,,\,t_0\,)\;=\;1\;.\;$ By lemma \ref{l1} we have, $\nu\,(\;x_0\;,\;\left\| {\;x\;} \right\|_{\,\alpha }^{\,\ast}\;)\;=\;1\,-\,\alpha\,$ for all $\alpha\,\in(\,0\,,\,1\,).\\$ Since $\nu\,(\,x_0\,,\,t_0\,)\;=\;1\;>\;1\,-\,\alpha$  it follows that $\nu\,(\;x_0\;,\;\left\| {\;x\;} \right\|_{\,\alpha }^{\,\ast}\;)\,\leq\,1\,-\,\alpha\,<\,\nu\,(\,x_0\,,\,t_0\,)\;$ and since $\nu\,(\,x_0\,,\,\cdot\,)$ is strictly non increasing $t_0\,<\;\left\| {\;x_0\;} \right\|_{\,\alpha }^{\,\ast}\;,\;\forall\;\alpha\,\in(\,0\,,\,1\,)\,.$ So $\nu\;^\prime\,(\,x_0\,,\,t_0\,)\;=\;\bigwedge\,\{\,1\,-\,\alpha\;:\;\left\| {\;x_0\;} \right\|_{\,\alpha }^{\,\ast}\;\leq\;t_0\}\;=\;1.$\\
Thus $\;\nu\,(\,x_0\,,\,t_0\,)\;=\;\nu\,^\prime(\,x_0\,,\,t_0\,)\;=\;1\,.\;$\\
\textbf{case-4:}\;$\;x_0\neq\theta,\,t_0\,(\,>0\,)\in \mathbb{R}$  such that $\nu\,(\,x_0\,,\,t_0\,)\;=\;0\;.\;$ From $(ix)$ it follows that, $\left\| {\;x_0\;} \right\|_{\,\alpha }^{\,\ast}<t_0,\;\;\forall\,\,\alpha\;\in(\,0,1\,)$. Therefore, $\left\| {\;x_0\;} \right\|_{\,\alpha }^{\,\ast}\;<\;t_0\;\;\Rightarrow\;\nu\,^\prime(\,x_0\,,\,t_0\,)=0,\;\;$  by (x).\\
Thus,  $\;\nu\,(\,x_0\,,\,t_0\,)\;=\;\nu\,^\prime(\,x_0\,,\,t_0\,)\;=\;0\,.$\\
$\textbf{case-5:} \;\;x_0\;\neq\;\theta\;,t_0\;(\,>\;0\,)\;\in\;\mathbb{R}$  such that $0\;<\;\nu\,(\,x_0\,,\,t_0\,)\;<\;1\,.$
Let, $\nu\,(\,x_0\,,\,t_0\,)\;=\;1\,-\,\beta\;,\;$ then from $(ix)$ we have
\begin{equation}\label{eq9}
\left\| {\;x\;} \right\|_{\,\beta }^{\,\ast}\;\leq\;t_0.
\end{equation}
Using (\ref{eq9}) from $(x)$ we get, $\nu\;^\prime\,(\,x_0\,,\,t_0\,)\;\leq\;1\,-\,\beta\,.\;$ Therefore,
\begin{equation}\label{eq10}
\;\nu\,(\,x_0\,,\,t_0\,)\;\geq\;\nu\,^\prime(\,x_0\,,\,t_0\,).\;
\end{equation}
Now from lemma \ref{l1} we have  $\nu\,(\,x_0\,,\,t_0\,)\;=\;1\,-\,\beta\;\Leftrightarrow\;\left\| {\;x\;} \right\|_{\,\beta }^{\,\ast}\;=\;t_0.\\$
Now, for $\;\beta\;<\;\alpha\;<\;1\;,\;$ let $\left\| {\;x\;} \right\|_{\,\alpha }^{\,\ast}\;=\;t\;^\prime\,.\;$ Then again by lemma \ref{l1} we have $\nu\,(\,x_0\,,\,t\,^\prime\,)\;=\;1\,-\,\alpha\,.\;$ So, $\nu\,(\,x_0\,,\,t\,^\prime\,)\;=\;1\,-\,\alpha\;<\;1\,-\,\beta\;=\;\nu\,(\,x_0\,,\,t_0\,)\,.\;$ Since, $\nu\,(\,x_0\,,\,\cdot\,)$ is strictly monotonically decreasing  and $\;\nu\,(\,x_0\,,\,t\,^\prime\,)\;<\; \nu\,(\,x_0\,,\,t_0\,)$ Therefore $\;t\,^\prime\;>\;t_0\,.\;$ Then for $\;\beta\;<\;\alpha\;<\;1\;,\;$ we have $\left\| {\;x\;} \right\|_{\,\alpha}^{\,\ast}\;=\;t\,^\prime\;>\;t_0\,.\;$ So,
\begin{equation}\label{eq11}
\nu\,^\prime(\,x_0\,,\,t_0\,)\;\geq\;1\,-\,\beta\;=\;\nu\,(\,x_0\,,\,t_0\,).
\end{equation}
Thus from (\ref{eq10}) and (\ref{eq11}) we have $\;\nu\,(\,x_0\,,\,t_0\,)\;=\;\nu\,^\prime(\,x_0\,,\,t_0\,)\,.\;$\\
Since, $\;(\,x_0\,,\,t_0\,)\;\in\;V\,\times\;\mathbb{R}$ is arbitrary $\;\nu\,^\prime(\,x\,,\,t\,)\;=\;\nu\,(\,x\,,\,t_0\,)\;\;,\;$ for all $(\,x\,,\,t\,)\;\in\;V\,\times\;\mathbb{R}\,.\\$

\begin{lemma}\label{l11}
In a fuzzy antinormed linear space $(\,V\,,\,A^*\,)$ with respect to a t-conorm $\diamond$ satisfying $(vi)$, $(vii)$ and $(viii)$, every sequence is convergent if and only if it is convergent with respect to its corresponding $\alpha$-norms, $\alpha\;\in\;(\,0\;,\;1\,).$
\end{lemma}
{\bf Proof.}$\;\;$
$\Rightarrow part:$ Let $(\,V\,,\,A^*\,)$ be a fuzzy antinormed linear space satisfying $(vi)$ and $(vii)$ and $\{x_n\}_{n\in\mathbb{N}}$ be a sequence in $V$ such that $x_n\;\rightarrow\;x$.
\[\mathop {\lim }\limits_{n\; \to \;\infty }\,\nu\,(\,x_n\,-\,x\,,\,t\,)\;=0\;\;,\;\forall\;t\,>\,0.\]
Choose $0\,<\,\alpha\,<\,1$. So, $\mathop {\lim }\limits_{n\; \to \;\infty }\,\nu\,(\,x_n\,-\,x\,,\,t\,)\;=0\;<\;1\,-\,\alpha\;\;\Rightarrow\;$ there exist $n_o\,(t)$ such that
\begin{equation}
\nu\,(\,x_n\,-\,x\,,\,t\,)\;<\;1\,-\,\alpha\;\;\forall\;n\,\geq\,n_0\,(t,\alpha).
\end{equation}
Now, \[\left\| {x_n\,-\,x} \right\|_{\,\alpha }^{\,\ast}=\bigwedge\,\{\,t>0:\;\nu\,(\,x_n\,-\,x\;,\;t\,)\;\leq\;1\,-\,\alpha\}\;
\Rightarrow\;\left\| {\;x_n\,-\,x\;} \right\|_{\,\alpha }^{\,\ast}\;\leq\;t\;,\forall\;n\,\geq\,n_0\,(t,\alpha).\]  Since, $t\,>\,0$ is arbitrary, \[\left\| {\;x_n\,-\,x\;} \right\|_{\,\alpha }^{\,\ast}\;\rightarrow\;0\;\;\; as\;\; n\rightarrow\,\infty\;\;,\;\forall\;\alpha\;\in\;(\,0\;,\;1\,).\]
$\Leftarrow part:$ Next, we suppose that, $\left\| {\;x_n\,-\,x\;} \right\|_{\,\alpha }^{\,\ast}\;\rightarrow\;0\;\;\; as\;\; n\rightarrow\,\infty\;\;,\;\forall\;\alpha\;\in\;(\,0\;,\;1\,)$.\\
Then for $\alpha\;\in\;(\,0\;,\;1\,)\;,\;\epsilon\,>\,0$ there exist $n_0\,(\alpha,\epsilon)$ such that
\begin{equation}
\left\| {\;x_n\,-\,x\;} \right\|_{\,\alpha }^{\,\ast}\;<\;\epsilon\;\;,\;\forall\;n\,\geq\,n_0\,(\alpha,\epsilon)\;\,,\;\alpha\;\in\;(\,0\;,\;1\,).
\end{equation}
Now,\[\nu\,(\,x_n\,-\,x\,,\,\epsilon\,)\;=\;\bigwedge\,\{\,1\,-\,\alpha\;:\;\left\| {\;x_n\,-\,x\;} \right\|_{\,\alpha }^{\,\ast}\;\leq\;\epsilon\}\hspace{8 cm}\]
\[\Rightarrow\;\nu\,(\,x_n\,-\,x\,,\,\epsilon\,)\;\leq\;1\,-\,\alpha\;\;,\;\forall\;n\,\geq\,n_0\,(\alpha,\epsilon)\;
\,,\;\alpha\;\in\;(\,0\;,\;1\,)\]
\[\Rightarrow\;\mathop {\lim }\limits_{n\; \to \;\infty }\,\nu\,(\,x_n\,-\,x\,,\,\epsilon\,)\;=\;0\hspace{6 cm}.\] Thus $x_n$ converges to $x$.

\begin{corollary}
Let $(V,A^\ast)$ be a fuzzy antinormed linear space with respect to a t-conorm $\diamond$ satisfying $(vi)$, $(vii)$ and $(viii)$. $W(\subseteq V)$ is closed in $(V,A^\ast)$ if and only if it is closed with respect to its corresponding $\alpha$-norms, $\alpha\in(0,1)$.
\end{corollary}

In the following lemma, a finite dimensional space is characterized by compact set in fuzzy environment and this will lead us to one of the fundamental differences between finite dimensional and infinite dimensional normed spaces with respect to fuzzy antinorms.

\begin{lemma}\label{rl}
(Riesz): Let $W$ be a closed and proper subspace of a fizzy antinormed linear space $(\,V\;,\;\nu\,)$ with respect to a t-conorm $\diamond$ satisfying $(vi)$, $(vii)$ and $(viii)$. Then for each $ \epsilon\;>\;0$ there exist $y\;\in\;V\,-\,W$ such that $\nu\,(\,y\;,\;1\,)\;\leq\;1\,-\,\alpha$ and $\nu\,(\,y\,-\,w\;,\;1-\epsilon\,)\;\leq\;1\,-\,\alpha$ for all $\alpha\;\in\;(\,0\;,\;1\,)$ and $w\;\in\;W$.
\end{lemma}
{\bf Proof.}$\;\;$
Recall that, $\left\| {\;x\;} \right\|_{\,\alpha }^{\,\ast}\;=\;\bigwedge\,\{\,t\;:\;\nu\,(\,x\;,\;t\,)\;\leq\;1\,-\,\alpha\}\;,\,\alpha\;\in\;(\,0\;,\;1\,)$ and $\{\,\left\| {\;.\;} \right\|_{\,\alpha }^{\,\ast}\;:\alpha\;\in\;(\,0\;,\;1\,)\,\}$ is an increasing family of $\alpha$-norms on a linear space $V$. Now by applying Riesz lemma for normed linear space, it follows that for any
$\epsilon\,>\,0$ there exist $y\,\in\,V\,-\,W$ such that
\begin{equation}
\left\| {\;y\;} \right\|_{\,\alpha }^{\,\ast}\;=\;1.
\end{equation}
\begin{equation}
\left\| {\;y\,-\,w\;} \right\|_{\,\alpha }^{\,\ast}\;>\;1-\epsilon\;\;\forall\;w\,\in\,W.
\end{equation}
Now, from theorem \ref{th1} for all $\alpha\;\in\;(\,0\;,\;1\,)$ we have
\[\nu\,(\,y\,,\,t\,)\;=\;\bigwedge\,\{\,1\,-\,\alpha\;:\;\left\| {\;y\;} \right\|_{\,\alpha }^{\,\ast}\;\leq\;t\}\hspace{11 cm}\]
\[\Rightarrow\;\nu\,(\,y\,,\,1\,)\;=\;\bigwedge\,\{\,1\,-\,\alpha\;:\;\left\| {\;y\;} \right\|_{\,\alpha }^{\,\ast}\;\leq\;1\}\;\hspace{2 cm}\]
\[\Rightarrow\;\nu\,(\,y\,,\,1\,)\;\leq\;1\,-\,\alpha.\;\hspace{5.7 cm}\]
Again,\[\nu\,(\,y\,-\,w\,,\,t\,)\;=\;\bigwedge\,\{\,1\,-\,\alpha\;:\;\left\| {\;y\,-\,w\;} \right\|_{\,\alpha }^{\,\ast}\;\leq\;t\}\hspace{4 cm}\]
\[\Rightarrow\;\nu\,(\,y\,-\,w\,,\,\epsilon\,)\;=\;\bigwedge\,\{\,1\,-\,\alpha\;:\;\left\| {\;y\,-\,w\;} \right\|_{\,\alpha }^{\,\ast}\;\leq\;\epsilon\}\;\hspace{2 cm}\]
\[\Rightarrow\;\nu\,(\,y\,-\,w\,,\,\epsilon\,)\;\leq\;1\,-\,\alpha.\;\hspace{6.5 cm}\]
Hence the proof.

\begin{theorem}
Let $(\,V\,,\,A^*\,)$ be a fuzzy antinormed linear space with respect to a t-conorm $\diamond$ satisfying $(vi)$, $(vii)$ and $(viii)$. If the set $\{\,x\;:\;\nu\,(\,x\,,\,1\,)\;\leq\;1\,-\,\alpha\}\;,\;\alpha\;\in\;(\,0\;,\;1\,)$ is compact then $V$ is a space of finite dimension.
\end{theorem}
{\bf Proof.}
It can be easily verified that $\{\,x\;:\;\nu\,(\,x\,,\,1\,)\;\leq\;1\,-\,\alpha\}\;\;
=\{\,x\;:\;\left\| {\;x\;} \right\|_{\,\alpha }^{\,\ast}\;\leq\;1\}\;,\;\alpha\;\in\;(\,0\;,\;1\,).$ By applying lemma \ref{rl}, it can be proved that if for some $\alpha\;\in\;(\,0\;,\;1\,)$ the set $\;\{\,x\;:\;\left\| {\;x\;} \right\|_{\,\alpha }^{\,\ast}\;\leq\;1\}\;$ is compact then $V$ is of finite dimensional. Using lemma \ref{l11}, it follows that, for some $\alpha\;\in\;(\,0\;,\;1\,)$, $\{\,x\;:\;\nu\,(\,x\,,\,1\,)\;\leq\;1\,-\,\alpha\}$ is compact then $V$ is a space of a finite dimensional.

%============================================
\section{Fuzzy $\alpha$-anti-convergence}
%============================================

In this section, the relations of fuzzy $\alpha$-anti-convergence, fuzzy $\alpha$-anti-Cauchyness, fuzzy $\alpha$-anti-compactness with respect to their corresponding increasing family norms are studied.\\

\begin{theorem}
Let  $(\,V\,,\,A^*\,)$ be a fuzzy antinormed linear space with respect to a t-conorm $\diamond$ satisfying $(vi)$, $(vii)$, $(viii)$ and $\{\;\left\| {\;\cdot\;} \right\|_{\,\alpha }^{\,\ast}\;:\;\alpha\;\in\;(\,0\;,\;1\,)\;\}$ be increasing family of norms of V, defined by $(ix)$. Then for any increasing (or, decreasing) sequence $\{\alpha_n\}_{n\in\mathbb{N}}$ in $(\,0\;,\;1\,),\;\;\alpha_n\;\rightarrow\;\alpha$ in $(\,0\;,\;1\,)$ implies $\;\left\| {\;x\;} \right\|_{\,\alpha_n }^{\,\ast}\;\rightarrow\;\left\| {\;x\;} \right\|_{\,\alpha }^{\,\ast}\;,\;\forall\;x\in V.$
\end{theorem}
{\bf Proof.}
For $x\;=\;\theta$, it is clear that $\alpha_n$ converges to $\alpha\;\;
\Rightarrow\;\left\| {\;x\;} \right\|_{\,\alpha_n }^{\,\ast}\;\rightarrow\;\left\| {\;x\;} \right\|_{\,\alpha }^{\,\ast}\,.$\\
Suppose $x\;\neq\;\theta$. Then from lemma \ref{l1}, for $x\neq \theta,\;\alpha\in(\,0,1\,)$ and $t\,^\prime\,>\,0,$ we have
\[\left\| {\;x\;} \right\|_{\,\alpha }^{\,\ast}\;=\;t\,^\prime\;\Leftrightarrow\;\nu\,(\,x_0\,,\,t\,^\prime\,)\;=\;1\,-\,\alpha.\]
Let $\{\alpha_n\}_{n\in\mathbb{N}}$ be an increasing sequence in $(\,0\,,\,1\,)$ such that $\alpha_n$ converges to $\alpha$ in $(\,0\,,\,1\,)$.\\
Let $\left\| {\;x\;} \right\|_{\,\alpha_n }^{\,\ast}=s_n\;$ and $\;\left\| {\;x\;} \right\|_{\,\alpha }^{\,\ast}=s.$\\
Then,
\begin{equation}\label{e41}
\nu\,(\;x\,,\,s_n\;)\;=\;1\,-\,\alpha_n\;\; and  \;\;\nu\,(\;x\,,\,s\;)\;=\;1\,-\,\alpha.
\end{equation}
Since $\{\;\left\| {\;\cdot\;} \right\|_{\,\alpha }^{\,\ast}\;:\;\alpha\;\in\;(\,0\;,\;1\,)\;\}$ is increasing family of norms and $\{s_n\}_{n\in\mathbb{N}}$ is increasing sequence of real numbers. Also, since $\{s_n\}_{n\in\mathbb{N}}$ is increasing sequence of real numbers and  bounded above by $s$, $\{s_n\}_{n\in\mathbb{N}}$ is convergent. Thus,
\begin{equation}\label{e42}
\mathop {\lim }\limits_{n\; \to \;\infty }\,\nu\,(\;x\,,\,s_n\;)\;=\;1\,-\mathop {\lim }\limits_{n\; \to \;\infty }\,\alpha_n\;\;\Rightarrow\;\nu\,(\;x\,,\,\mathop {\lim }\limits_{n\; \to \;\infty }\,s_n\;)\;=\;1\,-\,\alpha\;
\end{equation}
From (\ref{e41}) and (\ref{e42}) we have, $\;\nu\,(\;x\,,\,\mathop {\lim }\limits_{n\; \to \;\infty }\,s_n\;)\;=\;\nu\,(\;x\,,\,s\;)\;$ this implies $\mathop {\lim }\limits_{n\; \to \;\infty }\,s_n\;=\;s\;\;,\;$ by $(vii).\;\;$Therefore,\[\;\mathop {\lim }\limits_{n\; \to \;\infty }\,\left\| {\;x\;} \right\|_{\,\alpha_n }^{\,\ast}\;=\;\left\| {\;x\;} \right\|_{\,\alpha }^{\,\ast}\;\]
Similarly, if $\{\alpha_n\}_{n\in\mathbb{N}}$ is a decreasing sequence in $(\,0\,,\,1\,)$ and $\alpha_n$ converges to $\alpha$ in $(\,0\,,\,1\,)$ then it can be easily shown that $\;\left\| {\;x\;} \right\|_{\,\alpha_n }^{\,\ast}\;\rightarrow\;\left\| {\;x\;} \right\|_{\,\alpha }^{\,\ast}\;,\;\forall\;x\,\in\;V.$

\begin{definition}
Let $(\,V\;,\;A^\ast\,)$ be a fuzzy antinormed linear space with respect to a t-conorm $\diamond$ and $\alpha\;\in\;(\,0\;,\;1\,)$. A sequence $\{x_n\}_{n\in\mathbb{N}}$ in $V$ is said to be fuzzy $\alpha$-anti-convergent in $(\,V\;,\;A^\ast\,)$ if there exist $x\;\in\;V$ such that for all $t\,>\,0$,
\[\mathop {\lim }\limits_{n\; \to \;\infty }\,\nu\,(\,x_n\,-\,x \;,\; t\,) \;\,<\;\,1\,-\,\alpha.\;\;\;\;\] $x$ is called fuzzy $\alpha$-antilimit of $x_n$.
\end{definition}

\begin{theorem}
Let $(\,V\;,\;A^\ast\,)$ be a fuzzy antinormed linear space with respect to a t-conorm $\diamond$ satisfying $(vi)$ and $(viii)$. Then fuzzy $\alpha$-antilimit of a fuzzy $\alpha$-anti-convergent sequence is unique.
\end{theorem}
{\bf Proof.}
Let $\;\{x_n\}_{n\in\mathbb{N}}\;$ be a fuzzy $\alpha$-anti-convergent sequence and suppose it converges to $x$ and $y$ in $V$. Then for all $t\,>\,0\;,$
\[\mathop {\lim }\limits_{n\; \to \;\infty }\nu\,(\,x_n\,-\,x \;,\; t\,) \;\,<\;\,1\,-\,\alpha\;\;\;\;and \;\;\;\;\mathop {\lim }\limits_{n\; \to \;\infty }\nu\,(\,x_n\,-\,y\;,\; t\,) \;\,<\;\,1\,-\,\alpha.\]
Now,\[\nu\,(\,x\,-\,y\;,\;t\,)\;=\;\nu\,(\,x\,-\,x_n\,+\,x_n\,-\,y\;,\;t\,)\;\;,\;\forall\;n.\]
\[=\;\nu\,(\,x_n\,-\,x\;,\;t\,)\;\diamond\;\nu\,(\,x_n\,-\,y\;,\;t\,)\;\;,\;\forall\;n.\hspace{-4 cm}\]
Taking limit we have
\[\nu\,(\,x\,-\,y\;,\;t\,)\;=\;\mathop {\lim }\limits_{n\; \to \;\infty }\,\nu\,(\,x_n\,-\,x\;,\;t\,)\;\diamond\;\mathop {\lim }\limits_{n\; \to \;\infty }\,\nu\,(\,x_n\,-\,y\;,\;t\,)\]
\[\;<\;(\,1\,-\,\alpha\,)\;\diamond\;(\,1\,-\,\alpha\,)\;=\;(\,1\,-\,\alpha\,),\;\;(by\;(viii)).\]
That is, $\;\nu\,(\,x\,-\,y\;,\;t\,)<1,\;\;\;\;\forall\;t>0.\;.\;\;$\\ Therefore, $\;x\,-\,y\;=\;\theta\;\;\;\;\;by\;(vi).\\\;\Rightarrow\;x\;=\;y\;.$

\begin{theorem}\label{th41}
Let $(\,V\;,\;A^\ast\,)$ be a fuzzy antinormed linear space with respect to a t-conorm $\diamond$ satisfying $(vi)$ and $(viii)$. If $\{x_n\}_{n\in\mathbb{N}}$ be a fuzzy $\alpha$-anti-convergent sequence in $(\,V\;,\;A^\ast\,)$ such that $x_n$ converges to $x$. Then $\left\| {\;x_n\,-\,x\;} \right\|_{\,\alpha }^{\,\ast}\;\rightarrow\;0$ as $n\,\rightarrow\,\infty$.
\end{theorem}
{\bf Proof.}
Since $\;\{x_n\}_{n\in\mathbb{N}}\;$ be a fuzzy $\alpha$-anti-convergent sequence, suppose it converges to $x$, then for all $\;t\;>\;0\,,\\$
$\mathop {\lim }\limits_{n\; \to \;\infty }\,\nu\,(\,x_n\,-\,x\;,\;t\,)\;<\;\,1\,-\,\alpha.\\\;\Rightarrow\;\exists\;n_0(t)\;>\;0\;$ such that $\,\nu\,(\,x_n\,-\,x\;,\;t\,)\;<\;\,1\,-\,\alpha\;\;,\;\forall\;n\;\geq\;n_0(t).\\\;\Rightarrow\;\exists\;n_0(t)\;>\;0\;$ such that $\;\left\| {\;x_n\,-\,x\;} \right\|_{\,\alpha }^{\,\ast}\;\leq\;t\;\;\;\forall\;n\;\geq\;n_0(t).$\\Since $\;t\;>\;0\;$ is arbitrary, $\;\left\| {\;x_n\,-\,x\;} \right\|_{\,\alpha }^{\,\ast}\;\rightarrow\;0\;\;\; as \;n\;\rightarrow\;\infty\;.$

\begin{definition}
Let $(\,V\;,\;A^\ast\,)$ be a fuzzy anti-normed linear space with respect to a t-conorm $\diamond$ and $\alpha\;\in\;(\,0\;,\;1\,)$. A sequence $\{x_n\}_{n\in\mathbb{N}}$ in $V$ is said to be fuzzy $\alpha$-anti-cauchy sequence if \[\mathop {\lim }\limits_{n\; \to \;\infty }\;\nu\,(\,x_n\,-\,x_{n+p} \;,\; t\,) \;\,\leq\;\,1\,-\,\alpha\;\;\;\;\forall\,t\,>\,0\;,p=1,2,3,\cdots\]
\end{definition}
\smallskip
\begin{theorem}
Let $(\,V\;,\;A^\ast\,)$ be a fuzzy antinormed linear space with respect to a t-conorm $\diamond$ satisfying $(viii)$ and $\alpha\;\in\;(\,0\;,\;1\,)$. Then every fuzzy $\alpha$-anti-convergent sequence in $(\,V\;,\;A^\ast\,)$ is a fuzzy $\alpha$-anti-Cauchy sequence in $(\,V\;,\;A^\ast\,)$.
\end{theorem}
{\bf Proof.}
Let $\;\{x_n\}_{n\in\mathbb{N}}\;$ be a fuzzy $\alpha$-anti-convergent sequence and it converging to $x$. Then\[\mathop {\lim }\limits_{n\; \to \;\infty }\,\nu\,(\,x_n\,-\,x \;,\; t\,) \;\,<\;\,1\,-\,\alpha\;\]
Now, \[\nu\,(\,x_n\,-\,x_{n+p}\;,\;t\,)\;=\;\nu\,(\,x_n\,-\,x\,+\,x\,-\,x_{n+p}\;,\;t\,)\;\;,\;for\;\,p\;=\;1,2,3,\cdots.\]
\[\;=\;\nu\,(\,x_n\,-\,x\;,\;\frac{t}{2}\,)\;\diamond\;\nu\,(\,x_{n+p}\,-\,x\;,\;\frac{t}{2}\,)\;\;,\;for\;
\,p\;=\;1,2,3,\cdots.\hspace{-2 cm}\]
Therefore, \[\mathop {\lim }\limits_{n\; \to \;\infty }\;\nu\,(\,x_n\,-\,x_{n+p} \;,\; t\,) \;\,\leq\;\mathop {\lim }\limits_{n\; \to \;\infty }\;\nu\,(\,x_n\,-\,x\;,\;\frac{t}{2}\,)\;\diamond\;\mathop {\lim }\limits_{n\; \to \;\infty }\;\nu\,(\,x_{n+p}\,-\,x\;,\;\frac{t}{2}\,)\]
\[\;<\;(\,1\,-\,\alpha\,)\;\diamond\;(\,1\,-\,\alpha\,)\;=\;(\,1\,-\,\alpha\,),\;\;(by\;(viii)).\hspace{-0.7 cm}\]
Hence $\{x_n\}_{n\in\mathbb{N}}$ is a fuzzy $\alpha$-anti-Cauchy sequence in $(\;V\;,\;A^*\;)\,.$

\begin{note}
Every constant sequence in a fuzzy antinormed linear space $(\,V\;,\;A^\ast\,)$ with respect to a t-conorm $\diamond$ is a fuzzy $\alpha$-anti-Cauchy sequence in $(\,V\;,\;A^\ast\,)\;,\,\alpha\;\in\;(\,0\;,\;1\,)$.
\end{note}
{\bf Proof.}
Obvious.

\begin{theorem}
Let $(\,V\;,\;A^\ast\,)$ be a fuzzy antinormed linear space with respect to a t-conorm $\diamond$ satisfying $(vi)$ and $(viii)$. Then every Cauchy sequence in $(\;V\;,\;\left\| {\;\cdot\;} \right\|_{\,\alpha }^{\,\ast}\;)$ is a fuzzy $\alpha$-anti-Cauchy sequence in $(\,V\;,\;A^\ast\,)$, where $\;\left\| {\;\cdot\;} \right\|_{\,\alpha }^{\,\ast}\;$ denotes the increasing family of norms on $V\;$ defined by (ix), $\;\alpha\,\in\,(\,0\,,\,1\,).$
\end{theorem}
{\bf Proof.}$\;\;$
Choose $\;\alpha_0\;\in\;(\,0\;,\;1\,)\;$ arbitrary but fixed. Let, $\{y_n\}_{n\in\mathbb{N}}$ be a Cauchy sequence in $V$ with respect to $\;\left\| {\;\cdot\;} \right\|_{\,\alpha_0 }^{\,\ast}.\;$ Then \[\mathop {\lim }\limits_{n\; \to \;\infty }\;\left\| {\;y_n\,-\,y_{n+p}\;} \right\|_{\,\alpha_0 }^{\,\ast}\;=\;0\,.\] Then for any given $\,\epsilon\,(\,>\,0\,)\;$ there exist a positive integer $n_0(\epsilon)$ such that \\
$\;\left\| {\;y_n\,-\,y_{n+p}\;} \right\|_{\,\alpha_0 }^{\,\ast}\;<\;\epsilon\;\;,\;\forall\;n\;\geq\;n_0(\epsilon)$ and $p\;=\;1,2,3,\cdots.\\ \;\Rightarrow\;\bigwedge\,\{\;t\,>\,0\;:\;\nu\,(\;y_n\,-\,y_{n+p}\;,\;t\;)\;\leq\;1\,-\,\alpha_0\;\}\;<\;\epsilon\,. \\\;\Rightarrow\;$ there exist $t(\,n\,,\,p\,,\,\epsilon\,)\,<\,\epsilon\;$ such that\\ $\nu\,(\,y_n\,-\,y_{n+p}\;,\;t(\,n\,,\,p\,,\,\epsilon\,)\,)\,\leq\,1\,-\,\alpha_0\;\;,\;\forall\;n\;\geq\;n_0(\epsilon)$ and $p\;=\;1,2,3,\cdots.$
$\\\;\Rightarrow\;\nu\,(\;y_n\,-\,y_{n+p}\;,\;\epsilon\;)\;\leq\;1\,-\,\alpha_0\;$\\\\
Since $\;\epsilon\,(\,>0\,)$ is arbitrary,\\
$\mathop {\lim }\limits_{n\; \to \;\infty }\;\nu\,(\,y_n\,-\,y_{n+p} \;,\; t\,) \;\,\leq\;\,1\,-\,\alpha_0\;\;,\;\forall\;t\,>\,0\\ \;\Rightarrow\;\{y_n\}_{n\in\mathbb{N}}$  is fuzzy $\alpha_0$-anti-Cauchy sequence in $\;(\,V\;,\;A^\ast\,)\,.$\\
Since, $\alpha_0\;\in\;(\,0\,,\,1\,)$ is arbitrary, every Cauchy sequence in $\;(\,V\;,\;\left\| {\;\cdot\;} \right\|_{\,\alpha }^{\,\ast})\;$ is fuzzy $\alpha$-anti-Cauchy sequence in $\;(\,V\;,\;A^\ast\,)\;$ for each $\alpha\;\in\;(\,0\,,\,1\,)\,.$

\begin{definition}
Let $(\,V\;,\;A^\ast\,)$ be a fuzzy antinormed linear space with respect to a t-conorm $\diamond$ and $\alpha\;\in\;(\,0\;,\;1\,)$. It is said to be fuzzy $\alpha$-anti-complete if every fuzzy $\alpha$-anti-Cauchy sequence in $V$ fuzzy $\alpha$-anti-converges to a point of $V$.
\end{definition}

\begin{theorem}
Let $(\,V\;,\;A^\ast\,)$ be a fuzzy antinormed linear space with respect to a t-conorm $\diamond$ satisfying $(vi)$ and $(viii)$. If $(\,V\;,\;A^\ast\,)$ be fuzzy $\alpha$-anti-complete then $V$ be complete with respect to $\;\left\| {\;\cdot\;} \right\|_{\,\alpha }^{\,\ast}\;,\;\alpha\;\in\;(\,0\;,\;1\,).$
\end{theorem}
{\bf Proof.}$\;\;$
Choose $\;\alpha_0\;\in\;(\,0\;,\;1\,)\;$ arbitrary but fixed. Let $\{y_n\}_{n\in\mathbb{N}}$ be a Cauchy sequence in $V$ with respect to $\;\left\| {\;\cdot\;} \right\|_{\,\alpha_0 }^{\,\ast},\;$ then $\{y_n\}_{n\in\mathbb{N}}$ is fuzzy $\alpha_0$-anti-Cauchy sequence in $\;(\,V\;,\;A^\ast\,)\,.\;$\\
Since, $\;(\,V\;,\;A^\ast\,)\;$ is fuzzy $\alpha_0$-anti-complete, there exist $y\,\in\;V\;$ such that \\
$\mathop {\lim }\limits_{n\; \to \;\infty }\,\nu\,(\,y_n\,-\,y \;,\; t\,) \;\,<\;\,1\,-\,\alpha_0\;\;,\;\forall\;t\,>\,0\\\;\Rightarrow\;\mathop {\lim }\limits_{n\; \to \;\infty }\;\left\| {\;y_n\,-\,y\;} \right\|_{\,\alpha_0 }^{\,\ast}\;=\;0\;\;,\;$ by theorem \ref{th41}.\\
$\;\Rightarrow\;y_n\;\rightarrow\;y\;$ with respect to  $\;\left\| {\;\cdot\;} \right\|_{\,\alpha_0 }^{\,\ast}\;$.\\
$\;\Rightarrow\;(\;V\;,\;\left\| {\;\cdot\;} \right\|_{\,\alpha_0 }^{\,\ast}\;)\;$ is complete.\\
Since, $\alpha_0$ is arbitrary,  $(\;V\;,\;\left\| {\;\cdot\;} \right\|_{\,\alpha }^{\,\ast}\;)\;$ is complete.
\[\]{\bf Acknowledgments.\\}
The authors are grateful to the referees and the Editors for their fruitful comments, valuable suggestions and careful corrections in improving the paper in present form.

\end{document}